\documentclass{amsproc}

\usepackage {amsmath, amssymb, epsfig, enumerate}
\usepackage{a4wide}
\usepackage[latin1]{inputenc}
\newtheorem{thm}{Theorem}[section]
\newtheorem{prop}[thm]{Proposition}
\newtheorem{lem}[thm]{Lemma}
\newtheorem{cor}[thm]{Corollary}
\newtheorem{question}[thm]{Question}
\theoremstyle{definition}
\newtheorem{defi}[thm]{Definition}
\newtheorem{exam}[thm]{Example}
\theoremstyle{remark}
\newtheorem{rmk}[thm]{Remark}

\newcommand{\cc}{\mathbb{C}}
\newcommand{\bb}{\mathbb{B}}

\newcommand{\dd}{\mathbb{D}}

\newcommand{\nn}{\mathbb{N}}
\newcommand{\e}{\varepsilon}
\newcommand{\cv}{\rightarrow}

\newcommand{\fr}{\partial}
\newcommand{\om}{\Omega}
\newcommand{\set}[1]{\left\{#1\right\}}
\newcommand{\norm}[1]{\left\Vert#1\right\Vert}
\newcommand{\abs}[1]{\left\vert#1\right\vert}
\newcommand{\finc}{\subset \subset}
\newcommand{\cd}{\cc^2}
\newcommand{\pd}{{\mathbb{P}^2}}

\newcommand{\rest}[1]{ \arrowvert_{#1}}

\newcommand{\unsur}[1]{\frac{1}{#1}}
\newcommand{\el}{\mathcal{L}}

\newcommand{\tel}{{}^t \mathcal{L}}

\newcommand{\hl}{H{\'e}non-like }

\DeclareMathOperator{\supp}{Supp}
\DeclareMathOperator{\vol}{Vol}
\DeclareMathOperator{\Int}{Int}
\DeclareMathOperator{\jac}{Jac}

\DeclareMathOperator{\diam}{diam}

\begin{document}

\title[Connectivity of Julia sets in dimension 2]{Some remarks on the
  connectivity of Julia sets for
  2-dimensional diffeomorphisms.}

\author{Romain Dujardin}

\address{UFR de math{\'e}matiques,
Universit{\'e} Paris 7,
Case 7012,
2 place Jussieu,
75251 Paris cedex 05, France}
\email{dujardin@math.jussieu.fr}
\thanks{Part of this material was prepared as the contents of a lecture given
at a meeting of 
 the ACI Jeunes Chercheurs "Dynamique des applications polynomiales" in Rennes,  June 
2004. }
\date{\today}

\begin{abstract}
We explore the connected/disconnected  dichotomy for the Julia set
of polynomial
automorphisms of $\mathbb{C}^2$. We develop several aspects of the
question, which was first studied by Bedford-Smillie \cite{bs6,bs7}. 
We introduce a new sufficient condition for the connectivity of the Julia
set, that carries over for certain H{\'e}non-like and birational maps. 
We study the structure of disconnected Julia sets and the
associated invariant currents. This provides a simple approach to some results of
Bedford-Smillie, as well as some
 new corollaries --the connectedness locus is
closed, construction of external rays in the general case, etc. 

We also prove the following theorem: a hyperbolic polynomial
diffeomorphism of $\mathbb{C}^2$ with connected Julia set must have
attracting or repelling orbits. This is an analogue of a well known
result in one dimensional dynamics.
\end{abstract}
\maketitle

\section*{Introduction} \label{sec:intro}

If $p:\cc\cv\cc$ is a polynomial map, there exists a well known necessary and 
sufficient condition for the Julia set $J_p$ (or equivalently the
filled in Julia set 
$K_p$) to be connected. Namely $J_p$ is connected iff all critical points have bounded 
orbits. In case $p$ is hyperbolic, this is equivalent to saying that all critical points are 
attracted to periodic sinks (see \cite{cg} for a general account).\\

If now $f:\cd\cv\cd$ is a polynomial diffeomorphism, the Julia set has several
analogues. Let $K^+$ (resp. $K^-$) be the set of points with bounded forward 
(resp. backward) orbit. Let $J^\pm=\fr K^\pm$. It is known that $J^+$ 
is the Julia set of $f$ in the usual sense, for forward iteration. 
Let also $J=J^+\cap J^-$ and $J^*$ 
be the closure of the set of saddle orbits ($J^*\subset J$). \\

The sets $J^\pm$ are always connected. Indeed there exist invariant 
(1,1) currents
$T^\pm$ with $\supp( T^\pm)=J^\pm$, and with the additional property of
being extremal as positive closed currents.  
Extremality easily implies that the sets
$J^\pm=\supp (T^\pm)$ are connected. On the other hand, the sets $J^\pm$ are not locally connected in general, even if $f$ is hyperbolic on $J$. Indeed $J^\pm$ have
laminar structure and, for instance if $f(z,w)=(aw+z^2+c, az)$, with small $a$ and large 
$c$, $J^+$ has totally disconnected transversals. Later we shall see that the
connectivity properties 
of these transversals has to do with the connectivity of $J$ itself.\\

Eric Bedford and John Smillie gave in \cite{bs6} a necessary and sufficient 
condition for $J$ to be connected. Replacing $f$ by $f^{-1}$ if necessary 
we may assume the (constant) Jacobian determinant $\abs{{\rm Jac}\ f}$
is not larger than 1. 

\begin{defi}\label{def:uc}
$f$ is unstably connected if for some saddle point $p$, $W^u(p)\cap K^+$
has no compact component (for the topology induced by the isomorphism
$W^u(p)\simeq \cc$).
\end{defi}

An important issue in \cite{bs6} is that this condition is actually independent
of $p$.

\begin{thm}[\cite{bs6}]\label{thm:bs6}
If $f$ is a polynomial diffeomorphism with $\abs{{\rm Jac}\  f}\leq 1$, then 
$$ J \text{ is connected } \Leftrightarrow f \text{ is unstably connected}.$$ 
\end{thm}

A corollary, which was not realized there is that the connectedness locus
is closed in parameter space (corollary \ref{cor:closed}).

There is also a notion of unstable critical point, which plays the
 role of  escaping critical  points
 in one dimensional dynamics. Another result in \cite{bs6} is
$$ J \text{ is connected } \Leftrightarrow f \text{ has no unstable critical points}.$$
 
 A salient feature in that paper is that if $f$ is unstably connected, 
then $J^-\setminus
 K^+$ has the structure of a Riemann surface lamination. This allows to define external
 rays and study how landing of these influence the structure of $J$. \\
 
 In this article, we give a sufficient condition for the connectedness of the Julia set 
 of a polynomial diffeomorphism of $\cd$. This condition is not directly related
 to unstable connectedness and was considered in the context of hyperbolic maps in 
 \cite{bs7}. 
 
 We say $K^+$ is {\sl transversely connected} if there exists a holomorphic
 disk $V$, with $\fr V \cap K^+=\emptyset$ and $V\cap K^+\neq \emptyset$ (a
 transversal to $K^+$) such that $V\cap K^+$ is connected. We prove (\S \ref{sec:conn})
 that transversal
 connectedness implies $J$ is connected. One interest of this result is that it also
 provides 
 a simpler approach to the lamination structure of $J^-\setminus K^+$
 and clarifies somehow 
the analogy between  unstable and escaping critical points.
 Also the result is valid for certain
 H{\'e}non-like and birational maps. 
Notice that images of $\cc$  cannot be used in the context of
 H{\'e}non-like maps, so the approach of \cite{bs6} has to be modified.\\
 
 In \S\ref{sec:disc} we study the structure of the set $J^-$ for an unstably disconnected
 polynomial diffeomorphism. We prove that the unstable current $T^-$ is an integral
 of {\sl global} 
 submanifolds in a large bidisk $\bb$ containing $K$. This is strictly stronger 
 than being laminar but allows some folding. We can thus define external rays 
 in this case. We prove almost every ray lands and the landing measure 
 is the maximal entropy measure. 
 
 Another corollary (Proposition \ref{prop:extr} and corollary
\ref{cor:notextr}) is that  if $\bb$ is a sufficiently large bidisk,
 $$ f \text{ is unstably connected } \Leftrightarrow T^-\rest{\bb} \text{ is an
 extremal current.}$$

We  study the expansion/contraction of the Poincar{\'e} metric along
the leaves in this case, and derive some corollaries. 
This paragraph may also be seen as a gentle introduction 
to the notion of {\sl
  quasi-expansion} \cite{bs8}.
 
 In \S \ref{sec:hyp} we turn to hyperbolic polynomial diffeomorphisms. It is known 
 that if $f$ is hyperbolic on $J$ then ${\rm Int} K^+$ is the union of finitely many
 sink basins. We prove that if $f$ is unstably connected (or equivalently $J$ is 
 connected and $\abs{{\rm Jac} f}\leq 1$) then ${\rm Int} (K^+)$ is non empty, i.e. $f$ 
 does have 
 attracting periodic orbits. This  is analogous to the one dimensional
 case. This provides an alternate proof of the following fact
 \cite[Corollary A.3]{bs7}:
 the Julia set of a conservative and hyperbolic
 polynomial diffeomorphism cannot be connected.


\section{Maps with connected Julia sets}\label{sec:conn}

\subsection{Preliminaries on H{\'e}non-like mappings}\label{subs:prel_hl}
 Our treatment
of connectivity for polynomial diffeomorphisms was partly motivated 
by extending the results of \cite{bs6} to the \hl context. In particular
we will use the \hl formalism, even when considering H{\'e}non maps.

We begin with
some  notation: let  
$\bb$ be a bounded open set in $\cd$ biholomorphic to
a bidisk. We fix an isomorphism with the unit bidisk, and let
$\fr_v\bb$ (resp. $\fr_h \bb$) is the ``vertical''
(resp. ``horizontal'') part of the boundary, 
  $\fr_v\bb= \set{\abs{z}=1, \abs{w}<1}$ (resp.  $\fr_h
\bb= \set{\abs{z}<1, \abs{w}=1}$). 

\begin{defi}\label{def:hl}
A \hl mapping in $\bb$
is an injective holomorphic map $f:N(\overline\bb) \cv\cd$
(where $N(\overline\bb)$ is a neighborhood of $\overline\bb$)
such that $f(\bb)\cap \bb\neq\emptyset$, and satisfying:
\begin{itemize}
\item[\em i.] $f(\fr_v\bb)\cap \overline\bb = \emptyset$;
\item[\em ii.] $f(\overline \bb)\cap \fr \bb \subset \fr_v\bb$.
\end{itemize}
\end{defi}
 
Notice that the source and target bidisks may differ; this was the
case in the original definition of \cite{ho}.
A fundamental example of H{\'e}non-like map is  any polynomial automorphism
of $\cd$, while considered in a suitable large bidisk $\bb$. 
Indeed it is known that there exists a system of
coordinates $(z,w)$ so that $f$ is a composition of H{\'e}non maps 
$f_j(z,w)= (aw + p_j(z), z)$, and for every large enough $R$, $f$ is H{\'e}non-like in the 
bidisk $\bb=D(0, R)^2$. 

Perturbing such a map provides
many new examples. For instance consider the following
 polynomial birational map of $\cd$:
 $f(z,w)=((a+b(z))w +p(z), az)$ where $\abs{a}\leq 1$, and $b(z)$ is a
polynomial with  
${\rm deg}(b(z))\leq {\rm deg}(p(z))-1$. If $\norm{b(z)}_\bb$
 is small enough, $f$ is a small perturbation
of the H{\'e}non map $(aw+p(z), az)$, so $f$ is \hl in $\bb$
and it can be proved that  nonwandering dynamics occurs only in $\bb$. \\

We call {\sl horizontal} an object (form, current, subvariety) whose support 
stays away from the horizontal boundary $\fr^h\bb$. 
Vertical currents, as well as horizontal and vertical 
submanifolds, are defined 
analogously. The {\sl degree} of a horizontal subvariety $V$
is by definition the number of intersection points of $V$ with  a generic vertical line.
By definition, if $L$ is a horizontal line, $d={\rm deg}(f(L))$ is the degree of $f$.  
More generally if $V$ is a horizontal submanifold of degree $\deg(V)$,  
$f(V)\cap\bb$ is a horizontal submanifold of degree $d \deg(V)$.

There is  a corresponding 
notion of degree for a horizontal positive closed currents $T$. 
The degree is defined as the 
mass of the intersection (wedge product) of $T$ with a generic vertical line.
A {\sl normalized} current is a current of degree 1.

We now list some dynamical properties of \hl mappings \cite{hl}:
\begin{description}
\item[Invariant currents] $f$ acts by push forward  
on horizontal positive closed  currents. Let $\el=1_\bb \unsur{d}f_*$
be the associated graph transform operator. If $T$ is any horizontal normalized
positive closed current, the sequence $\el^nT$ converges to the unique normalized
$\el$ invariant current $T^-$. Moreover $T^-$ has laminar structure and continuous 
potential.
Similar results hold for pull backs of vertical 
currents.
\item[Invariant measure] $\mu=T^+\wedge T^-$ is the unique measure of maximal
entropy $\log d$. It is mixing, hyperbolic, and describes the asymptotic distribution
of periodic orbits.
\end{description}

If we let $K^\pm=\set{x\in\bb,~\forall n\geq 0~f^{\pm n}(x)\in
\bb}=\bigcap_{n\geq 0}f^{\mp n}(\bb) $ and $J^\pm=\fr K^\pm$, an interesting 
 open question is  whether $\supp(T^\pm)=J^\pm$. 
Equality holds for
polynomial diffeomorphisms of $\cd$.

For the birational perturbations of H{\'e}non maps considered
above, the equality $J^+=\supp (T^+)$ is true. This is an easy consequence of
the existence of the 
 rate of escape function $g^+$, which is psh, nonnegative, continuous, and such
that $dd^c g^+=T^+$ and $K^+=\set{g^+=0}$. Here $K^+$ is the set of points with 
bounded orbits, and agrees in $\bb$ with the previously defined  $K^+$. \\

\subsection{A result on the intersection of positive closed currents}
We begin with a result on the support of wedge product of positive closed currents, 
which has some independent interest. It is a very flexible
generalization of \cite[Proposition 2.3]{bs3} and \cite[Lemma 5.4]{bs6}.

\begin{thm}\label{thm:support}
Let $\om\subset \cd$ be a bounded  open set, and $T_1=dd^cu_1$, 
$T_2=dd^cu_2$ be positive closed currents
in $\om$ such that the wedge product $T_1\wedge T_2$ is admissible. 
We moreover assume 
\begin{equation}\label{eq:hypsupport}
\emptyset\neq \supp (T_1)\cap \supp(T_2)\finc \om.
\end{equation} Then 
 $\displaystyle{\int_\om T_1\wedge T_2>0}$. 
\end{thm}

Remark that if $\om$ is exhausted by pseudoconvex open sets, hypothesis
(\ref{eq:hypsupport}) implies that the wedge product $T_1\wedge T_2$
is well defined \cite{sib}. The theorem has the following corollary.

\begin{cor} \label{cor:support}
If $f$ is a \hl map, $\supp(\mu)$ intersects every connected 
component of $\supp(T^+)\cap\supp(T^-)$.
\end{cor} 

\begin{proof} We may assume $\om$ is connected. 
Let $\chi$ be a test function with 
$\chi=1$ in a neighborhood of
 $\supp(T_1)\cap\supp(T_2)$ and $\supp(d\chi)\cap\supp(T_1)\cap
\supp(T_2)=\emptyset$. 
Then 
\begin{equation}\label{eq:mass}
\int_\om T_1\wedge T_2= \int \chi dd^cu_1\wedge dd^cu_2= \int
u_1 dd^c\chi\wedge dd^cu_2 <\infty
\end{equation}
since $u_1$ is pluriharmonic near
$\supp(dd^c\chi)\cap\supp(T_2)$. Let $T_1^\e=dd^cu_1^\e$ be the standard 
regularization of $T_1$. $u^\e_1$ is obtained by the convolution of $u_1$ with
a radial approximation of $\delta_0$, with support in $B(0,\e)$.
Reducing $\om$ slightly if necessary we may assume 
$T_1^\e$ is well defined on $\om$. 

Now $u_1-u_1^\e=0$ in the set $\set{p\in \om, d(p, \supp(T_1))>\e}$. In particular if 
$\e$ is small enough,$u_1-u_1^\e=0$ on  $\supp(dd^c\chi)\cap\supp(T_2)$
and by (\ref{eq:mass}) we get that $\int T_1\wedge T_2=\int T_1^\e\wedge T_2$.
We do the same for $T_2$, hence  $\int T_1\wedge T_2=\int T_1^\e\wedge T_2^\e$.\\

The next step is to prove that for small $\e>0$, $\supp (T_1^\e)\cap \supp(T_2^\e)\neq
\emptyset$ (of course $\supp (T_1^\e)\cap \supp(T_2^\e)\finc\om$). 
Once again it is enough to show that  $\supp (T_1^\e)\cap \supp(T_2)\neq
\emptyset$. First, it is classical that 
there exists a nonnegative
plurisubharmonic function $v$ in $\om\backslash \supp(T_2)$
tending to $+\infty$ on $\supp(T_2)$. Indeed for every $\e>0$, the function 
$u_2^\e-u_2$ is psh, nonnegative on $\om\setminus\supp(T_2)$, and 
positive in the neighborhood of $\supp(T_2)$. Pick a decreasing
sequence $\e_k\cv 0$ and just take 
$$v= \sum_{k=0}^\infty C_k(u^{\e_k}_2-u_2),$$ where  the constants $C_k$ are adjusted 
so that $v\cv+\infty$ on $\supp(T_2)$ (the sum is locally finite in 
$\om\backslash \supp(T_2)$).\\

Assume next  that for every $\e>0$, $\supp (T_1^\e)\cap \supp(T_2)=\emptyset$. We will
use a Kontinu{\"a}tzsatz-type argument. There exists a constant $c$ such that
in a neighborhood $N$ of $\fr\om$, $\supp(T_1)\cap N\subset \set{v\leq c }$. This 
also holds for $T_1^\e$ for small $\e$, with a fixed $c$. Now if 
$\supp (T_1^\e)\cap \supp(T_2)=\emptyset$, $v$ is well defined on $\supp (T_1^\e)$
and by the maximum principle (see below), 
 $\supp(T_1^\e)\subset \set{v\leq c }$. Letting
$\e\cv0$, we get a contradiction. 

The  maximum principle on a positive closed current is classical. The
proof goes as follows: assume $w$ is a smooth strictly
psh function in $U$, such that 
$w < 0$ on $\fr U\cap \supp (T)$, and $w>0$ at some point $p\in
\supp(T)\cap U$. Then $w^+=\max(w,0)$ vanishes near $\supp(T)\cap \fr
U$ and $w^+=w$ near $p$. Then  
$$0< \int_U dd^c w^+\wedge T= \int w^+\wedge dd^cT =0;$$ a contradiction.\\ 

So by now we may assume that we are in the conditions of the theorem, with 
smooth $T_1$ and $T_2$. Smoothness is not enough to ensure 
$\int T_1\wedge T_2>0$, since $T_1$ and $T_2$ might be only semipositive 
on  $\supp(T_1)\cap\supp(T_2)$. Assume this is so,
and pick $p\in \supp(T_1)\cap\supp(T_2)$. 
Consider a small connected neighborhood $N$
of $id$ in the group $SU_2$ of rotations around $p$, and the current
$$T'_1= \int_N r_* T_1 dr$$ where $dr$ is the normalized Haar measure 
in $N$. Elementary linear algebra shows  that $T'_1$ is strictly positive at $p$. 
so $\int T'_1\wedge T_2>0$. 

It  remains to prove that 
$\int T'_1\wedge T_2=\int T_1\wedge T_2$. For every rotation
$r$ in $N$, $r$ is connected to the identity by a path in $N$, so
$r_*T_1$ is homotopic to $T_1$. The homotopy formula (e.g. \cite{sim}) 
asserts that 
$$r_*T_1-T_1= d h_*(T_1\otimes[0,1]) + h_*(dT_1\otimes [0,1]),$$ 
where $h(\cdot, t)$ is the homotopy connecting $id$ and $r$, and 
$T_1=T_1\mathbf{1}_\om$ 
is viewed as a current with boundary in the neighborhood of $\overline\om$.
 We infer
$$\int \chi(r_*T_1-T_1)\wedge T_2 = \int \chi d h_*(T\otimes[0,1]) \wedge T_2
+ \int  \chi h_*(dT_1\otimes [0,1]) \wedge T_2.$$ 
The first integral on the right hand side vanishes because 
$d\chi = 0$ near the support of $\left(h_*(T_1\otimes[0,1]) \wedge T_2\right)$
and the second because 
$\chi=0$ near $\supp\left(h_*(dT_1\otimes [0,1])\right) $.
\end{proof}


\subsection{Polynomial automorphisms} We begin the discussion 
on connectivity by considering 
transverse connectivity for polynomial automorphisms 
of $\cd$. Fix a polynomial 
diffeomorphism $f$ of degree $d$ in $\cd$ with non trivial dynamics. 
We will prove transversal connectedness of $K^+$
implies $J$ is connected.  We also directly recover 
the lamination structure of $J^-\backslash K^+$, and its ergodicity. 

\begin{defi} \label{def:tc}
A transversal to $K^+$ is  a holomorphic disk $V$ such that 
$\fr V \cap K^+=\emptyset$ and $V\cap K^+\neq \emptyset$. $K^+$ is said to be 
transversely connected if there exists a transversal $V$ such that
 $K^+\cap V$ is connected.
\end{defi}

Notice that if $K^+\cap V$ has an isolated connected component (e.g. if $K^+\cap V$
has finitely many components) then $K^+$ is transversely connected. In \cite{bs7}
it was proved that if $f$ is uniformly hyperbolic then $K^+$ is 
transversely connected iff $J$ is 
connected. Here we extend the ``only if'' statement to the general case.  

\begin{thm}\label{thm:connex}
Let $f$ be a polynomial diffeomorphism of $\cd$ and assume $K^+$ is 
transversely connected. Then the Julia set $J$ is connected. 
Moreover $J^-\setminus K^+$ supports a unique Riemann surface lamination 
 which is uniquely  ergodic.  
\end{thm}

The proof will consist in several steps. Let $\bb=D(0,R)^2$ be a large bidisk; 
$f\rest{\bb}$ is a H{\'e}non-like map.  
We also introduce $V^+=\set{(z,w)\in \cd, \ \abs{z}>\abs{w}, \ \abs{z}>R}$.
It is classical that $f(\bb)\subset \overline{\bb}\cap V^+$ and $f(V^+)\subset V^+$.

\begin{proof} We fix a transversal $V$ such that $V\cap K^+$ is connected. Since $\fr V
\cap K^+=\emptyset$, $\fr V$ escapes under iteration. So there exists $n_0\geq 1$
such that $f^{n_0}(V)\cap\bb$ is a horizontal submanifold in $\bb$. Since 
$f^n(V\cap K^+)=f^n(V)\cap K^+$ is connected,  replacing $V$ by $f^n(V)\cap\bb$, we
may assume $V$ is a horizontal submanifold in $\bb$ (i.e. $\fr V\subset  \fr^v\bb$), 
of degree $\deg(V)$. Here ``submanifold'' means a complex submanifold without 
boundary.\\

\noindent{\bf Step 1.}  If $V$ is as above, then for every $n\geq 1$,
 $f^n(V)\cap \bb$ is connected.\\

Let $\bb_{-n}=\bb\cap f^{-1}(\bb)\cap\cdots \cap f^{-n}(\bb)$; from the fact that
$f(V^+)\subset V^+$ one deduces  $f^n(V)\cap\bb= f^n(V\cap\bb_{-n})$.
 It then suffices to prove that for every $n$, 
$V\cap \bb_{-n}$ is connected.
For this, we just remark that $K^+$ intersects every connected
component of 
$V\cap \bb_{-n}$. 
Indeed, if $U$ is such a component, then $U$ is a (boundaryless)
submanifold of $\bb_{-n}$. So $f^n(U)$ is a non trivial horizontal submanifold in
$\bb$. It then 
follows from the Stokes theorem and degree considerations
that $\int T^+\wedge [f^n(U)] >0$  (see e.g. \cite[Prop. 2.7]{hl} or
theorem \ref{thm:support}), so $f^n(U)$ intersects
$K^+$.\\

\noindent{\bf Step 2.} The laminar structure of $J^-\setminus K^+$.\\

Recall first the {\sl Riemann-Hurwitz formula}. Let $\Delta$ is a simply
connected horizontal 
submanifold in $\bb$ of degree $\delta$ (so $\Delta$ is a union of $k$ holomorphic
disks), 
 and $\pi$ denote the first projection $(z,w)\mapsto z$.  Then the number of critical 
 points of $\pi\rest{\Delta}$ 
 (vertical tangencies), counted with multiplicity, equals $\delta-k$. \\
  
 Consider the transversal $V$ as before; we have seen that 
 that for $n\geq 1$, $f^n(V)\cap\bb$ is a connected horizontal submanifold. 
 $f^n(V)\cap\bb=f^n(V\cap\bb_{-n})$ is a planar surface, and  the 
 maximum principle actually implies it is a disk. So by the Riemann-Hurwitz formula
the number of vertical tangencies  on $f^n(V)\cap\bb$ equals $d^n\deg(V)-1$.

Now increase $R$ to get a larger bidisk $\bb'$. For large enough $n$, 
$f^n(V)\cap\bb'$ is also a horizontal disk in $\bb'$, of the same degree 
$d^n\deg(V)$, and the same number of vertical tangencies. Thus all vertical tangencies
of $f^n(V)\cap\bb'$ are inside $\bb$, or equivalently, the projection
$$\pi: f^n(V)\cap(\bb'\setminus \bb) \longrightarrow D(0,R')\setminus D(0,R)$$ is 
a covering. Hence for every simply connected open subset $Q\subset 
D(0,R')\setminus D(0,R)$, $f^n(V)\cap\pi^{-1}(Q)$ is the union of $d^n\deg(V)$ 
graphs. Recall that $(d^n\deg(V))^{-1}[f^n(V)\cap\bb']$ converges to the unstable 
current $T^-$. It then follows (see e.g. \cite{bs5}) that 
$T^-\rest{\pi^{-1}(Q)}$ is a {\sl uniformly laminar current}, made up of integration
currents over the limiting graphs. 

We  have thus proved that in $\bb'\setminus\bb$,
$J^-=\supp T^-$ supports a lamination $\el^-$. Now if $\bb'$ is so large that 
$(\bb'\setminus\bb)\cap J^-$ contains a fundamental domain for the action 
of $f$ on $J^-\backslash K^+$, we obtain that $J^-\setminus K^+$ is laminated.
The uniqueness of the lamination is obvious since $J^-$ has empty interior.\\

\noindent{\bf Step 3.} Unique ergodicity of the transverse measure.\\

The lamination $\el^-$ constructed above carries a foliation cycle $T^-$, 
so it has an invariant transverse measure. We prove it is uniquely
ergodic, that is, the  positive invariant transverse measure is unique.
It implies of course ergodicity, 
i.e. any two transversals of positive measure are connected by holonomy. 
The lamination will actually be (uniquely) ergodic in each $\bb'\setminus\bb$. 

We use here \cite[Prop. 2.13]{bs6}, 
which is itself a slight extension of a result in \cite{fs}. The claim
is  that any positive
closed current $S$ supported in $K^-\setminus\bb$ is a multiple of $T^-$. In 
particular, $T^-\rest{\cd\setminus\bb}$ is extremal in $\cd\setminus\bb$. Since 
any invariant transverse measure induces a foliation cycle, which is a positive closed
current supported in $K^-\setminus\bb$, we get that the unique (up to a scalar multiple)
invariant transverse measure is the one induced by $T^-$.\\

Now if $\nu$ is an invariant transverse measure for the lamination
 $\el^-\rest{\bb'\setminus \bb}$, we prove it can be extended to a transverse measure 
on $\el^-\rest{\cd\setminus\bb}$, 
hence it is again induced by $T^-$. Recall that for every simply 
connected open subset $Q$ in $\cc\setminus D(0,R)$, the leaves of the lamination 
in $\pi^{-1}(Q)$ are graphs over $Q$. Hence the  vertical lines $\pi^{-1}(p)$, 
$p\in\cc\setminus D(0,R)$ are global transversals. Given any 
two points $p\in D(0,R')\setminus D(0,R)$ and $q\in \cc\setminus D(0,R)$, one may 
thus transport the transverse measure from $\pi^{-1}(p)$ to 
$\pi^{-1}(q)$ by using a simply connected $Q\subset  \cc\setminus D(0,R)$ containing 
$p$ and $q$. So the transverse measure $\nu$ extends from $\el^-\rest{\bb'\setminus 
\bb}$ to $\el^-\rest{\cd\setminus \bb}$ and we get the desired conclusion.\\

Remark that since $T^-$ has full support in $J^-$, the lamination is also minimal.\\

\noindent{\bf Step 4.} Connectivity of $J$.\\

%


By corollary \ref{cor:support}, every connected component 
of $J$ contains a point of
$J^*=\supp(\mu)$, where $\mu=T^+\wedge T^-$ is the maximal entropy measure.
Similarly to \cite{bs6}, the connectivity of $J$ will
follow from the following fact  ``for every $\e>0$, 
any two points $p$ and $q$ in $J^*$ are joined by a path lying in the 
$\e$-neighborhood $J_\e$ of $J$''. \\

Indeed,  fix  $p,q\in J^*$, and $\e>0$. 
For every $x\in J^-\setminus K^+$, $f^{-n}(x)$ converges to $J$, 
so by compactness  there exists
an integer $n$ such that $f^{-n}((J^-\cap \bb) \setminus K^+)\subset J_\e$. 
By Pesin Theory, $\mu$ almost  every
 point has a local unstable manifold $W^u_{loc}$,
subordinate to a piece of unstable lamination with positive measure.
Moving $p$ and $q$ slightly is necessary, we may assume this is true for $p$ and $q$.
 Moreover, since the laminar structure of $T^-$ is subordinate to
the decomposition in 
unstable manifolds \cite{bls}, 
for every local unstable manifold $W^u_{loc}$, 
$W^u_{loc}\cap (J^-\backslash K^+)$, which is non empty, is subordinate 
to the lamination $\el^-$ constructed above.

Let $M=\sup_\bb\norm{df}$.
Since  the lamination $\el^-$ is uniquely ergodic in 
a small neighborhood of $\fr^v\bb$, its restriction to  $\bb\setminus K^+$ is
ergodic, so we can find 
a path $\gamma$ subordinate to a leaf  of $\el^-\rest{\bb\setminus K^+}$,
 and joining two points
$p_1$ and  $q_1$, with  $d(f^n(p),p_1)<\frac{\e}{M^n}$ and. 
$d(f^n(q),q_1)<\frac{\e}{M^n}$. Hence $f^{-n}(\gamma)$ is a path contained in 
$J_\e$ and joining $f^{-n}(p_1)$ and $f^{-n}(q_1)$, which are respectively $\e$-close 
to $p$ and $q$.
\end{proof}

\begin{rmk} \label{rmk:crit}
There is a strong analogy between our condition and the condition of non escaping 
critical points in one variable dynamics. The role of ``critical point'' is played here by
vertical tangencies. We indeed proved in step 2 that vertical tangencies do not
escape a certain bidisk. See also the analysis in example \ref{ex:hyp} below. 

A more precise notion of escaping critical point is developed in \cite{bs5,bs6}. 
Escaping critical points are the critical points of $G^+$ restricted to unstable
manifolds. Using the fact that the critical points of $G^+\rest{V}$ are 
the vertical tangencies for the invariant projection $\varphi^+$ 
(the ``B{\"o}ttcher coordinate'') it is
not difficult to prove that if $V$ is a transversal to $K^+$
$$V\cap K^+ \text{ is connected }\Leftrightarrow G^+\rest{V} \text{ has no  
critical points}.$$
\end{rmk}

In the next proposition we relate transverse connectivity to unstable
connectivity. It does not follow  from Theorem \ref{thm:bs6} 
because of the assumption on the Jacobian.

\begin{prop}\label{prop:tcuc}
If $K^+$ is transversely connected, then $f$ is unstably connected.
\end{prop}

From \cite[Corollary 7.4]{bs6} one deduces:

\begin{cor}\label{cor:tcuc}
If $K^+$ is transversely connected, then $\abs{Jac(f)}\leq 1$.
\end{cor}

\begin{proof} Assume $f$ is unstably disconnected. So for some saddle periodic
point $p$, $W^u(p)\cap K^+$ has a compact component $C$. The action 
of $f$ on $W^u(p)\simeq\cc$ (where 0 stands for $p$) is a nontrivial dilatation. So $0\notin C$, and $\bigcup_{n\geq 1} f^{-n}(C)$ consists of infinitely many compact 
components of $W^u(p)\cap K^+$, close to 0.

On the other hand $[V]\wedge T^+>0$, so 
we claim that the stable manifold $W^s(p)$ has
a transverse intersection point with $V$. A proof goes as follows
(we use ideas from \cite[\S 9]{bls}): $\unsur{d^n}
f^n_*[V] \cv cT^-$, where $c=\int[V]\wedge T^+$. Consider  
a Pesin box $P$ of positive  measure; the local stable manifolds
associated to points form a piece of 
stable lamination of positive transverse measure, hence 
a uniformly laminar current  \cite[\S 8]{bls}
$S^+\leq T^+$, such that $T^-\wedge S^+>0$, 
because $G^-$ cannot be harmonic on $W^s_{loc}(p)$.
 It is classical that $S^+$ then has continuous potential. Thus 
$$\unsur{d^n} f^n_*[V]\wedge S^+\cv T^-\wedge S^+>0,$$ and we get that 
for large $n$, $f^n(V)$ has
intersection points with a set of positive measure of
disks in $S^+$ and most of them are transverse
\cite[Lemma 6.4]{bls}. 

Assume $f^n(V)$ intersects $W^s_{loc}(x)$, for some $x\in P$, not
necessarily periodic.
By Poincar{\'e} recurrence, we
 may suppose that for infinitely many $n_j$, $f^n_j(x)\in P$, so by
using the stable manifold theorem, we get that for  large $j$, 
$f^{n+n_j}(V)$ contains a disk arbitrarily close to $W^u_{loc}(f^{n_j}x)$. Now if $p$ is 
any saddle point, $W^s(p)$ has transverse intersection points with any set of positive
measure of unstable disks \cite[Lemma 9.1]{bls}. We conclude
that $W^s(p)$ must intersect $f^{n+n_j}(V)$, hence $V$, transversely.\\

For $N>>1$, consider a set $C_1, \ldots C_N$ of open and closed subsets 
of $W^u(p)\cap K^+$ close to $p$. For each $C_i$, consider a simple
curve $\gamma_i$ 
enclosing $C_i$, and such that $\gamma_i \cap K^+=\emptyset$. Choose 
local coordinates $(x,y)$ so that $W^u_{loc}(p)\subset (y=0)$, and  
extend $\gamma_i$
to a piece of 3-submanifold
 $\widetilde\gamma_i$, transverse to $W^u_{loc}(p)$,
 by adding vertical holomorphic disks. We assume the vertical disks 
 are so small that $\widetilde\gamma_i\cap K^+=\emptyset$

By the Lambda lemma, for large $n$, 
$f^n(V)$ contains a graph over $W^u_{loc}(p)$, and very close to it. So 
the curves $\widetilde\gamma_i\cap f^n(V)$ cut out $N$ 
disks $\Delta_1\ldots \Delta_N$ 
in $f^n(V)$. Of course as $n\cv\infty$ the disks $\Delta_i$ converge to 
${\rm Int}(\gamma_i)$. This forces $\Delta_i\cap K^+$ to be non empty 
for large $n$: for instance use the fact that $[{\rm Int}(\gamma_i)]\wedge T^+>0$ 
(see step 1 above)
and the continuity of the potential $G^+$.

We have proved that in any neighborhood of any transverse intersection point
of $W^s(p)\cap V$, $K^+\cap V$ has at least $N$ connected components, 
which clearly contradicts transverse connectivity.
\end{proof}

\begin{rmk}
The proof of the proposition provides another approach to the fact that unstable 
connectivity is independent of the chosen unstable manifold.
\end{rmk}

In the next proposition we  use the fact that  
if $f$ is unstably connected, $J^-\setminus K^+$  has the structure described 
in Theorem \ref{thm:connex}
above (a uniquely ergodic lamination with foliation cycle $T^-$). 

\begin{prop}\label{prop:extr}
If $f$ is unstably connected then $T^-\rest{\bb}$ is an extremal current.
\end{prop}

\begin{proof} We have seen that the unstable current $T^-$ is uniformly laminar
in $\cd\setminus K^+$, and moreover if $R$ and $R'$ are sufficiently large 
the transverse measure in $\el^-\rest{\pi^{-1} (A(R,R'))}$ is ergodic, where 
$A(R,R')$ denotes the annulus $D(0,R')\setminus D(0,R)$, and $\pi$ is the first 
projection.
Any positive closed current $S\leq T^-\rest{\pi^{-1} (A(R,R'))}$ 
is uniformly laminar and subordinate to $T^-$ (this is a result about
analysis on laminations, see \cite{l2}). By ergodicity, we
conclude that  $T^-\rest{\pi^{-1} (A(R,R'))}$ is extremal.

Let now $S$ be  a positive 
closed current in $\bb$, with $S\leq T$. In particular $S=cT^-$ ($c\leq 1$) in the 
neighborhood of $\fr^v\bb$. Consider the positive closed current 
$\widetilde{S}$ on $\cd$ defined by $  \widetilde{S}=S$ in $\bb$, and 
$\widetilde{S}=cT^-$ outside $\bb$. Since $T^-$ is extremal in $\cd$, 
$\widetilde{S}=c T^-$ everywhere, and we conclude that $T^-\rest{\bb}$ is extremal.
\end{proof}

\begin{exam}\label{ex:hyp} Among  well understood 
  polynomial diffeomorphisms are perturbations of one dimensional hyperbolic 
maps. Here we show that if $p$ is a hyperbolic polynomial with
connected Julia set, then for small $\abs{a}$, the H{\'e}non map 
$f_a(z,w)=(aw+p(z), az)$ is transversely connected. 
The proof is inspired  by  \cite{ho}, from which the result may
actually be extracted.  \\

Let $R$ be so large that  $\abs{z}=R$ implies $\abs{p(z)}>2R$, and let
$\bb =D(0,R)^2$. First, if $\abs{a}\leq 1$, $f$ is a H{\'e}non-like map of
degree $d$ in $\bb$ (see lemma \ref{lem:perturb} below).
We will show that if $L$ is any horizontal line in $\bb$, and $a$ is
small enough, $L\cap K^+$ is connected. Since for such a $L$
$$L\cap K^+=\bigcap L\cap \bb_{-n}=\bigcap \overline{L\cap\bb_{-n}},$$ 
and every component of $L\cap\bb_{-n}$ intersects $K^+$, this is
equivalent to saying that for every $n$, $L\cap \bb_{-n}$ is
connected. By the Riemann-Hurwitz formula (see step 2 above) this is
in turn equivalent to the fact that $f^n(L)\cap \bb$ has $d^n-1$
vertical tangencies: ``critical points do not escape $\bb$''.

We begin with a useful lemma.

\begin{lem}\label{lem:perturb}
Let $U_1$ and $U_2$ be two topological disks such that
$\overline{U_i}\subset D(0,R)$. Assume that 
\begin{itemize}
\item[-] $p:p^{-1}(U_2)\cv U_2$ is a branched cover of some degree
  $k$.
\item[-] ${\rm dist}(p(\fr U_1),U_2)>\delta$.
\end{itemize}
Then $f_a:U_1\times D(0,R)\cv U_2\times D(0,R)$ is a H{\'e}non-like
(crossed) map of degree $k$ as soon as $\abs{a}\leq
\min(\frac{\delta}{R},1)$.
\end{lem}
 
\begin{proof} Items {\em i.} and {\em ii.} of definition \ref{def:hl}
easily hold. We check that $f_a: U_1\times D(0,R)\cv U_2\times D(0,R)$ 
has degree $k$: if $u_2\in U_2$, we need to prove that
$f^{-1}(\set{u_2}\times D(0,R))$, which is a vertical submanifold
in $U_1\times D(0,R)$, has degree $k$. It is defined by the equation
$$\set{(z,w),\ aw+p(z)=u_2},$$ and clearly its intersection with 
$(w=0)$ has $k$ points counted with multiplicity.
\end{proof}

Assume now that $p$ is a hyperbolic polynomial with connected Julia 
set.  Fix a neighborhood $N$ of $J$ such that 
$p^{-1}(N)\finc N$, and $p:p^{-1}(N)\cv N$ is a covering of degree $d$.
$p$ is a strict expansion for the Poincar{\'e} metric of $N$. We denote by
$U(z,\e)$ the ball of radius $\e$ around $z$ for the Poincar{\'e} metric of $N$. 
Fix $\e$ small enough, so that for  every $z\in J$, $U(z,\e)$ is a topological 
disk, and $p$ is univalent on $U(z,\e)$. In particular,
$$p:p^{-1}(U(p(z),\e))\cap U(z,\e)\cv U(p(z), \e)$$ is a biholomorphism. 
Reducing $\e$ once again, one  may further assume that 
$${\rm dist}\big(p(\fr U(z,\e)), U(p(z), \e)\big)>\delta,$$ for some constant
$\delta$ independent of $z$.  By the preceding lemma, 
for $\abs{a}\leq \frac{\delta}{R}$, 
$$f_a: U(z,\e)\times D(0,R)\cv U(p(z),\e)\times D(0,R)$$ is  a
H{\'e}non-like map of degree 1. Let $U=\bigcup_{z\in J} U(z,\e)$.\\

Let us now consider a horizontal line $L$ in $\bb$. By the above argument,
for every $\abs{a}\leq \frac{\delta}{R}$, and every
$n\geq 1$, all iterates $f_a^n(L)$ are graphs over $U$: indeed 
since $p:p^{-1}(N)\cv N$ has degree $d$, the only contribution
to $f_a^n(L)$ over $U$ comes from $U\times D(0,R)$. Now for fixed
$n$, if $a\cv 0$, the vertical tangencies of $f_a^n(L)$ converge to 
the $d^n-1$ (with multiplicities) critical values of $p^n$, located in $K$. Since 
the vertical tangencies cannot cross $U\times D(0,R)$, we conclude that in $\bb$
$f_a^n(L)$ has $d^n-1$ vertical tangencies in $K\times D(0,R)$, hence 
$f_a^n(L)\cap \bb$ is connected. \hfill $\square$
\end{exam}

\begin{rmk} It can be proved that under these assumptions $J^+\cap \bb$ is 
laminated by vertical holomorphic disks, moving holomorphically with $a$. 
Moreover 
when $a\cv 0$, the lamination converges to the trivial lamination 
of $J\times D(0,R)$. From this one concludes that any slice 
$J^+(f_a)\cap L$  is the image of $J(p)$ by a holomorphic motion. 
\end{rmk}

\subsection{\hl mappings} We present some  connectedness results in
the H{\'e}non-like setting\footnote{The main idea in the 
proof of the next theorem
originates from a remark made to me by N. Sibony.}. 
The picture is much less precise than in
 the case of polynomial automorphisms, 
in particular we cannot prove that  points in $J$ are connected by paths
subordinate to unstable manifolds. 

Let $f$ be a \hl map in $\bb$. Notice first that if 
$f$ is dissipative or conservative (i.e. $\abs{{\rm Jac}(f)}\leq 1$), then
$K^-$ has measure zero, hence $J^-=K^-$. Indeed $f^{-1}K^-\finc K^-$,
so if $f$ contracts volumes $K^-$ cannot have positive measure.

If $f$ is a birational perturbation of a dissipative H{\'e}non map, as
considered in \S\ref{subs:prel_hl}, the equality $J^+=\supp T^+$ holds.

\begin{thm}\label{thm:connex_hl}
Let $f$ be a H{\'e}non-like map of degree $d>1$ in $\bb$, and assume that 
$J^-=K^-$.
\begin{enumerate} 
\item If  $J^+=\supp (T^+)$ and if $K^+$ is transversely connected, 
then $J=J^+\cap J^-$ is connected.
\item If $J^-=\supp(T^-)$ and if there exists a transversal $V$ to
  $\supp(T^+)$ such that $\supp T^+ \cap V$ is connected, then
$\supp(T^+)\cap\supp(T^-)$ is connected.
\end{enumerate}
\end{thm}

\begin{proof}{\em (1)} If $V$ is a transversal to $K^+$, then  for $n$
large enough, $\el^n[V]$ is a non trivial horizontal positive closed current in
$\bb$ (where $\el$ is the graph transform operator for currents, see
\S \ref{subs:prel_hl}). It follows that the sequence of currents
$(\el^n[V])$ converges to $cT^-$, with $c>0$. Let $\mu_n=\el^n[V]\wedge
T^+= (f^n)_*([V]\wedge T^+)$; since $T^+$ has continuous potential, $\mu_n\cv c\mu$ as
$n\cv\infty$. 

Since $K^+$ is holomorphically convex and $V$ is a holomorphic
disk, $V\setminus
K^+$ has no compact components. It follows that $J^+\cap V=\supp
T^+\cap V$ is connected, and so does $f^n(J^+\cap V)$. 

On the other hand if  $J_\e$ is the $\e$-neighborhood 
of $J$, for large $n$ one has
$$\supp \mu_n\subset \supp(T^+) \cap f^n(V) = f^n(J^+\cap V)\subset
J^+\cap K^-_\e\subset J_\e.$$
 We will conclude that $J=J^+\cap J^-$ is connected by 
showing that every connected component of $J^+\cap J^-$ contains a
point of $\supp(\mu)$.

So let $J_1$ be an open and closed subset of $J$, $J_1=J\cap \om$,
with $\fr\om\cap J=\emptyset$. We prove that $\mu(\om)>0$. Since
$J^-=K^-$,  the points in $J^+\cap
\fr\om=\supp(T^+\rest{\overline\om})\cap \fr\om$ escape under
backwards iteration. So if $\tel$ denotes the pull back graph
transform operator, $\tel^n(T^+\rest{\om})$ is closed and vertical
for large $n$, and non trivial since $\supp(\tel^n(T^+\rest{\om}))\cap
J^-\neq \emptyset$. In particular $\tel^n(T^+\rest{\om})\wedge T^->0$,
hence $(T^+\rest{\om})\wedge T^->0$. \\

{\em (2)} The reasoning is similar. Let $\psi$ be a nonnegative
test function 
in $V$, with $\psi=1$ near $\supp(T^+)$, so that $\psi[V]\wedge
T^+=[V]\wedge T^+$.
Assume for a moment that 
$\el^n\psi[V]\cv cT^-$, with $c>0$. Then $\mu_n =\el^n\psi[V]\wedge
T^+$ converges to $c\mu$, and since $\supp(T^-)=J^-=K^-$, for every 
$\e>0$
$$\supp \mu_n\subset \supp(T^+) \cap f^n(V) =f^n (\supp(T^+)\cap V) 
\subset (\supp (T^+)\cap \supp (T^-))_\e$$ for large $n$. We conclude 
by using corollary \ref{cor:support}. 

It remains to prove our claim. The difficulty is that we can not
assume $\el^n[V]$ is closed after a few iterations. Nevertheless 
$c=\int T^+\wedge [V]>0$ by theorem \ref{thm:support}, and 
it will be a consequence of \cite{hl} and \cite{dds} that
$\el^n(\psi[V])\cv cT^-$.  

An easy adaptation of proposition 4.8, and theorem 4.10   
in \cite{hl} shows that if the mass of the sequence of currents
$\el^n(\psi[V])$ is bounded, then $\el^n(\psi[V])\cv cT^-$.  
Let $\Theta$ be a smooth vertical positive closed current. 
From \cite[Proposition 4.13]{dds}, one may write
$\tel^n\Theta=dd^cu_n$, where $(u_n)$ is a uniformly bounded sequence 
of psh functions. Then the sequence
$$ \int \el^n(\psi[V])\wedge\Theta = \int_V \psi dd^cu_n = \int_Vu_n dd^c\psi$$
is uniformly bounded, and since in the neighborhood of $K^+$, masses
may be evaluated by using vertical closed positive currents we get the result.
\end{proof}

\begin{rmk} ~
\begin{itemize}
\item[-] If $f$ is a polynomial diffeomorphism, by corollary
  \ref{cor:tcuc}, the assumption $J^-=K^-$ is a consequence of
 transversal connectedness. 
\item[-] In case {\em (1)}, one may reproduce steps 1 and 2 in the
  proof of theorem \ref{thm:connex}, and infer that $J^-$ is a 
lamination near $\fr^v\bb$. 
\end{itemize}
\end{rmk}

\section{Maps with disconnected Julia sets}\label{sec:disc}

\subsection{The connectedness locus is closed}
Throughout this section we assume $f$ is a polynomial automorphism of
$\cd$ of degree $d$, and we fix as before a bidisk $\bb$ such that $f\rest{\bb}$
is  H{\'e}non-like  of degree $d$.

Recall that $f$ is said to be unstably disconnected if for some saddle
point $p$, $W^u(p)\cap K^+$ has a compact component (for the leafwise 
topology). We saw in
proposition \ref{prop:tcuc} that this is independent of the saddle
point $p$. Recall also that if $\abs{\jac(f)}\leq 1$, $J$ is
disconnected iff $f$ is unstably disconnected. Moreover, if 
$\abs{\jac (f)}<1$, $f$ is always stably disconnected, and if
$\abs{\jac(f)}=1$, $f$ is stably disconnected iff 
$f$ is unstably disconnected \cite[Corollary 7.4]{bs6}. 

\begin{lem}\label{lem:disk}
If $p$ is any saddle point, 
$W^u(p)\cap K^+$ has a compact component iff $W^u(p)\cap \bb$ has a
relatively compact component.
\end{lem}

\begin{proof} Without loss of generality assume $p$ is fixed.
Let first $C$ be a compact component of  $W^u(p)\cap K^+$, and
$\gamma$ be a Jordan curve in $W^u(p)\setminus K^+$
surrounding $C$, i.e. $C\subset \Int(\gamma)$. Then $\gamma$ escapes
under iteration, and for some $n$, $f^n(\Int(\gamma))\cap\bb$ is a horizontal
disk in $\bb$, relatively compact in $W^u(p)$. 

Conversely assume $C$ is a relatively compact component of $W^u(p)\cap
\bb$. Then $C$ is a non trivial horizontal submanifold in $\bb$.
Thus $C$ intersects $K^+$, and $C\cap K^+\finc C$ is a compact
component of $W^u(p)\cap K^+$.
\end{proof}

We obtain as a corollary that the connectedness locus is closed in
parameter space\footnote{This result stemmed out during 
a conversation with Eric Bedford and Misha Lyubich at the Snowbird
conference in holomorphic dynamics in June 2004.}. 
Because the  parameter space of polynomial
diffeomorphisms is not well understood, we state the result in terms
of 1-parameter families.

\begin{cor}\label{cor:closed}
Let $\set{f_\lambda, \lambda\in\Lambda}$ 
be a holomorphic 1-parameter family of polynomial diffeomorphisms. 
Then the connectedness locus
 $$\set{\lambda\in \Lambda, \ J(f_\lambda) \text{ is connected}}$$
is closed in $\Lambda$.
\end{cor}
 
\begin{proof} We prove the disconnectedness locus is open. So suppose
  that  at $0\in \Lambda$,  $J(f_{0})$ is disconnected. Reducing
  $\Lambda$ slightly if necessary, we may assume all $f_\lambda$ are
H{\'e}non-like in $\bb$.

$f_0$ is both stably and unstably disconnected. We prove stable
disconnectedness persists under perturbation.
By the preceding lemma, there is a saddle point $p_0$,
  such that $W^u(p_0)\cap \bb$ has a compact component. 
Saddle points and their unstable manifolds
 vary continuously (holomorphically) 
under perturbation, so for $\lambda$ close to $0$, there is
  a holomorphic family $p_\lambda$ of saddle points corresponding to
  $p_0$, and every relatively compact part of $W^u(p)$ can be followed 
holomorphically. So for $\lambda$ close to $0$,  $W^u(p_\lambda)\cap
  \bb$ has a compact component and we are done.
\end{proof}

\subsection{Structure of $T^-$} 
The currents  $T^\pm$ have laminar structure in general. We will 
 show that if $f$ is unstably disconnected, $T^-\rest{\bb}$ is an
 integral of horizontal submanifolds of $\bb$. In other words, almost
 every leaf of $T^-\rest{\bb}$ is a finite branched cover over the basis
$\dd$ for the natural vertical projection. 

Recall that $f$ is always
 stably or unstably disconnected, so the result always apply to at
 least one of $T^+$ or $T^-$.\\

If $T$ is a laminar current, we say that a disk $\Delta$ is
subordinate to $T$ if there exists a non trivial uniformly laminar current 
$S\leq T$ with $\Delta$ as a plaque. See \cite{structure} for
definitions and related results. We begin with a useful proposition.

\begin{prop}\label{prop:subord}
If $p$ is a saddle point, any relatively compact disk $\Delta\subset
W^u(p)$ is subordinate to $T^-$.
\end{prop}

\begin{proof} If $W^u_{loc}(p)$ denotes any neighborhood of $p$ in
  $W^u(p)$, pulling back $\Delta$ by $f$ if necessary, 
it is enough to prove  $W^u_{loc}(p)$ is subordinate to $T^-$.

We saw in the course of the proof of proposition \ref{prop:tcuc} that
if $Q$ is a  Pesin box, the local unstable lamination $\el^u(Q)$
supports a non trivial uniformly laminar current subordinate to
$T^-$. Moreover if $p$ is any saddle point, $W^s(p)$ has transverse
intersection points with
$\el^u(Q)$. Thus there exists a uniformly laminar current 
$$S_0= \int [\Delta_\alpha]d\nu_0(\alpha)\leq T^-$$ 
where the $\Delta_\alpha$ are 
 small disks intersecting $W^s(p)$ transversely at one point. 

By the Lambda lemma for every $\alpha$, $f^{n}(\Delta_\alpha)$
contains graphs over $W^u_{loc}(p)$, arbitrary close to $p$ when $n$
is large. So there exists a uniformly laminar current $S_j$
subordinate to $d^{-j}(f^j)_* S_0\leq T^-$, and  made up of graphs over 
$W^u_{loc}(p)$, close to $p$. 
Extracting a subsequence if necessary we may assume
the $S_j$ have disjoint supports. Define a uniformly laminar current
$S$ in the neighborhood of $p$ as  $S=\sum c_jS_j$, where the
constants $c_j\leq 1$ are adjusted so that the series converges.
$W^u_{loc}(p)$ is subordinate to $S$, and $S\leq T^-$ so this solves
the problem.
\end{proof}

The structure theorem for the current $T^-\rest{\bb}$ follows easily.
It is stated in $\bb$, however the bidisk may be arbitrary large.
A consequence is that 
$J^-$ is locally the limit of the union of an increasing
sequence of laminations. 
 
\begin{thm}\label{thm:structure}
If $f$ is unstably disconnected, then  
there exists a family of uniformly laminar currents $T^-_k$ in $\bb$, 
respectively made up of 
horizontal disks of degree $k$, such that $$T^-\rest{\bb}=\sum_{k=1}^\infty T_k^-.$$
\end{thm}

\begin{proof} Lemma \ref{lem:disk} asserts that for any unstable
  manifold $W^u(p)$, $W^u(p)\cap \bb$ has a relatively compact
  component. Such a component is a horizontal submanifold of some
  degree $k$ in $\bb$. By proposition \ref{prop:subord} this
  horizontal submanifold is subordinate to $T^-$, so there exists a
  uniformly laminar $T_0$ made up of submanifolds of degree $k$
in $\bb$, with $0\leq T_0\leq T^-\rest{\bb}$. Notice that 
by the maximum principle, all
  submanifolds involved here are holomorphic  disks.

The sequence of cut-off iterates $\unsur{d^n}(f^n)_*
(T_0)\rest{\bb}\leq T^-$
converges to $T^-\rest{\bb}$, and each of these currents is uniformly
laminar and made up of global submanifolds of bounded degree.
Any uniformly laminar current subordinate to $T$ is of the form 
$hT$, where $0\leq h\leq 1$ is a measurable function, 
locally constant along the leaves.
So we can  write  $\unsur{d^n}(f^n)_*
(T_0)\rest{\bb}$ as $h_n T^-\rest{\bb}$, where $0\leq h_n\leq 1$, and the convergence
result says that  $h_n\cv 1$, $T^-$-a.e.
Hence the  sequence of currents   
$$T_n = \max\left(T_0, \ldots, {d^{-n}}(f^n)_*(T_0)\rest{\bb}
\right)= \max(h_0, \ldots, h_n) T^-\rest{\bb}$$ increases to $T^-\rest{\bb}$. 
It is clear that each $T_n$ is uniformly laminar, and made of
submanifolds of bounded degree, 
so $T^-$ has the desired structure. 
\end{proof}

Since none of the $T^-_k$ can be extremal we get the following corollary.

\begin{cor}\label{cor:notextr}
If $f$ is unstably disconnected, then $T^-\rest{\bb}$ is not an
extremal current.
\end{cor}

\begin{rmk} The decomposition of $T^-$ into an integral of extremal
  components can easily be deduced from the structure theorem. It
  appears that the extremal components of $T^-\rest{\bb}$ are
  irreducible horizontal submanifolds.

In case $f(\bb)\cap\bb$ is disconnected, 
a similar corollary already appears in \cite[\S 5.2]{dds}. More specifically,
if $f(\bb)\cap\bb =U_1\cup\cdots\cup U_m$ ($m>1$)
each point in $K^+$  (resp. $K^-$)
 can be assigned an itinerary  $\alpha\in \set{1,\ldots, m}^\nn$. 
This induces decompositions of the currents $T^\pm\rest{\bb}$ in terms
of itinerary sequences.

 Each open  set $U_i$ comes equipped with a multiplicity $d_i=\deg f(L)\cap
 U_i$, where $L$ is a horizontal line in $\bb$. If at least one
 partial degree $d_i$ is larger than one, then, with respect to a
 natural measure on  $\set{1,\ldots, m}^\nn$ related to the $d_i$,
for almost every symbolic
 sequence $\alpha$, $d_{\alpha(0)}\cdots d_{\alpha(n)}\cv\infty$. From
 this it follows that in general the decomposition of $T^-\rest{\bb}$
 given in \cite{dds} is not the extremal decomposition of $T^-$.  
\end{rmk}

Another corollary is the construction of external rays in the unstably
disconnected case. External rays for unstably connected maps were
constructed in \cite{bs6} using the lamination structure of
$J^-\setminus K^+$. 

\begin{cor}\label{cor:externalrays}
There exists a measured family $(\mathcal{E}, \nu)$
of external rays, defined as gradient
lines of $G^+$ restricted to unstable manifolds. 
Almost every ray lands and if $e$ denotes the endpoint mapping, then $e_*\nu$ equals  the
maximal entropy  measure $\mu$.
\end{cor} 

\begin{proof} Pick up a generic leaf $M$ of
  $T^-\rest{\bb}$, so $M$ is a horizontal disk in $\bb$ of finite degree.
Reducing $M$ a little  if necessary, it  has finite volume, and we may
  assume $\fr M$ if of the form $\set{G^+=r}$ for some $r>0$, and does
  not intersect the critical set of $G^+\rest{M}$. 

On $M\setminus K^+$, $G^+\rest{M}$ is a positive harmonic function, so outside
the critical points of $G^+\rest{M}$, one may flow along the gradient
lines (in the sense of decreasing $G^+$). If $\set{G^+=s}$ is
conveniently oriented, the 1-form
$d^c(G^+\rest{M})\rest{\set{G^+=s}}$ defines a positive 
 measure on $M\cap \set{G^+=s}$ which is invariant by the flow. If
$s<r$, only finitely many gradient lines issued from $\set{G^+=r}$ hit
a critical point of $G^+\rest{M}$ before attaining $\set{G^+=s}$. This
 defines a measurable map $e_{r,s}:\set{G^+=r}\cv\set{G^+=s}$ such
 that
 $$(e_{r,s})_*d^c(G^+)\rest{\set{G^+=r}\cap
 M}=d^c(G^+)\rest{\set{G^+=s}\cap M}.$$

That almost every ray lands is a consequence of the fact that $M$ has
finite volume, see \cite[\S 7]{bj} for a proof.\\

 At the level of the whole unstable lamination, we get 
a map $e_{r,s}:\set{G^+=r}\cv\set{G^+=s}$, defined almost everywhere
with respect to the measure $d^cG^+\rest{\set{G^+=r}}\wedge dd^cG^-$. 
Indeed, notice that $d^cG^+\rest{\set{G^+=r}}\wedge dd^cG^-$ is the
integral of the family of measures  $d^cG^+\rest{\set{G^+=r}\cap M}$
with respect to the transverse measure induced by $T^-$.
One has 
 $$(e_{r,s})_*d^cG^+\rest{\set{G^+=r}}\wedge dd^cG^-
=d^cG^+\rest{\set{G^+=s}}\wedge dd^cG^-. $$ On the  
other hand \cite{bs6}
 $$d^cG^+\rest{\set{G^+=r}}\wedge dd^cG^-=dd^c\max(G^+, r)\wedge
 dd^cG^- \underset{r\cv 0}{\longrightarrow} \mu=dd^cG^+\wedge dd^cG^-, $$
 so the landing measure is the maximal entropy measure.
\end{proof}

\subsection{The Poincar{\'e} metric along the leaves} We are interested in
 the expansion properties of $f$ with respect to the Poincar{\'e} metric
 along the leaves of $T^-\rest{\bb}$. We will use the structure 
of $T^-\rest{\bb}$ to derive (rough) estimates on the Poincar{\'e} metric.

Notice that we cannot rule out the possibility of singular
leaves on a set of zero transverse measure (the disks constructed in
theorem \ref{thm:structure} are smooth). 
This does not affect the treatment of the Poincar{\'e} metric,
since if $\Gamma$ is a singular irreducible holomorphic disk of degree
$k$, there exists a parametrization $\phi:\dd\overset{\sim}\cv\Gamma$
such that $\pi\circ\phi$ is a branched cover of degree $k$. Indeed
just consider the uniformization $\dd\cv\widehat{\Gamma}$ of the
desingularization of $\Gamma$. Pushing forward by $\phi$ 
 defines unambiguously the Poincar{\'e} metric on $\Gamma$. Here we
shall treat horizontal holomorphic disks regardless of their regular or
 singular nature; this actually has to be considered only in
 proposition \ref{prop:nogrowth}.  \\
 
Assume $M$ is a horizontal disk of  degree $k$ in $\bb$. Let $\bb_0$
and $\bb_1$ be ``vertical sub-bidisks'' of the form
$\bb_i=\bb\cap \set{\abs{z}<R_i}$, with $R_0<R_1$. It is classical that there
exists a  constant $C$ such that 
$$k\pi  R_1^2 \leq \vol(M\cap \bb_1) \leq   Ck\pi R_1^2.$$
The following nice result is \cite[Theorem 3.1]{bs8}.

\begin{thm}
Let $M$ be a horizontal disk of degree $k$
as above. Assume $M_0 \subset M_1$ are respective 
 connected components of  $M\cap\bb_0$ and $M\cap \bb_1$. Then
 $M_1\setminus M_0$ is  an annulus and there exists a
 constant $C$ depending only on $R_0$ and $R_1$ such that 
$${\rm Modulus}(M_1\setminus M_0)\geq \frac{C}{k}$$
\end{thm}

We denote by $\rho_M$ the Poincar{\'e} (Kobayashi) metric of the
submanifold $M$ in $\bb$.

\begin{cor}\label{cor:Poincare}
There exist  constants $C$ and $\lambda>1$, depending on $(k,R_0,R_1)$ 
such that the hyperbolic diameter of $M_0$ in $M_1$
is bounded by $C$ and $\rho_{M_0}\geq \lambda \rho_{M_1} \geq \rho_M$.  
\end{cor}

\begin{proof} We may think of $M_1$ as being the unit disk $\dd$. If
  $E\finc \dd$ define \cite[\S 2.3]{mcm}
$${\rm Mod}(E,\dd)=\sup\set{{\rm Modulus}(A), \ A\subset \dd
\text{ is an annulus enclosing }E}.$$ Then the hyperbolic diameter of
$E$ and ${\rm Mod}(E,\dd)$ are related by 
$$(\diam_\dd(E)\cv\infty) \Leftrightarrow ( {\rm Mod}(E,\dd)\cv 0).$$
So by the previous theorem there exists a constant $C(k,R_0,R_1)$ such
that  $$\diam_\dd(M_0)\leq C.$$ Moving $M_0$ by an isometry of $\dd$ if
necessary, 
we may assume $0\in M_0$. The last estimate implies then the existence of
a radius $r(k,R_0,R_1)<1$ such that $M_0\subset D(0,r)$, and the
result follows.  \\

We may actually give an explicit estimate: when the diameter is large the following holds
$$C \unsur{\diam_\dd(E)}\geq {\rm Mod}(E,\dd) \geq C
e^{-\diam_\dd(E)},$$
so if $k$ is large $\diam_\dd(E)\leq C k$. Since 
$$\rho(0, z)=\unsur{2}\log\frac{1+\abs{z}}{1-\abs{z}}$$ 
  we get $r\leq 1-e^{-Ck}$, hence the estimate 
$$\rho_{M_0}\geq \unsur{1-e^{-Ck}}\ \rho_{M_1} \geq \rho_M$$ 
on the Poincar{\'e} metric.  
\end{proof}

Using non uniform expansion along the leaves for the Poincar{\'e} metric
allows to recover the following result from \cite{bs6}. Notice that
the current $T^+$ induces by restriction a measure  on every unstable
leaf. 
These measures are equivalent to the ``unstable conditionals'' of
$\mu$ \cite[Prop 3.1]{bls}. 

\begin{thm}\label{thm:totdisc}
If $f$ is unstably disconnected, then for $\mu$-a.e. $p$, almost every
component of $K^+\cap W^u(p)$ is a point, that is, for
$T^+\rest{W^u(p)}$-a.e. $x$, the connected component of $x$ in
$K^+\cap W^u(p)$ is $\set{x}$. 
\end{thm}

\begin{proof} For a $\mu$-generic point $p$, let $K^{+,u}(p)$ be the connected
component of  $K^+\cap W^u(p)$ containing $p$. Since the unstable
conditionals of $\mu$ are induced by $T^+$ it is enough  to prove that
for $\mu$-a.e. $p$,  $K^{+,u}(p)=\set{p}$. 

Theorem \ref{thm:structure} provides us with a
  decomposition $T^-\rest{\bb}=\sum T^-_k$. Fix an integer $K$ so that 
$T^-_{\leq K}=\sum_{k\leq K} T^-_k$ is nontrivial. The measure 
$\mu_{K}= T^-_{\leq K}\wedge T^+\leq \mu$ has positive mass, and there 
exists a set $E$ of total $\mu_K$ mass, such that if
$p\in E$, the connected component $\Gamma^u(p)$
of $W^u(p)\cap \bb$ containing $p$ is a horizontal disk of degree
$k\leq K$. since $K$ is arbitrary it suffices to
 show that for $p\in E$, $K^{+,u}(p)=\set{p}$.

There exists a radius $R_1$ such that $f^{-1}(\bb)\cap \bb \subset
\bb\cap \set{\abs{z}<R_1}$.
So by the previous corollary, if $p\in E$, the hyperbolic diameter of 
$K^+\cap \Gamma^u(p)$ is uniformly bounded by a constant $D$, and 
if $p$, $f^k(p)\in E$, $f^{-k}: \Gamma^u(f^k(p))\cv \Gamma^u(p)$ is a
contraction of factor a least $\lambda$. 
Moreover for almost every
$p\in E$ there are infinitely many $n$ such that $f^{n}(p)\in E$. Label
these values as $n_j$, $j\geq 1$.
For such a $n_j$, 
$f^{-n_j}:\Gamma^u(f^{n_j}(p))\cv \Gamma^u(f(p))$ is a  contraction of 
factor $\lambda^j$ for the Poincar{\'e} metric. Of course $f(K^{+,u}(p))=K^{+,u}(f(p))$.
Thus the hyperbolic diameter of  
$K^{+,u}(p)$ is bounded by $D/\lambda^j$ and the result
follows.
\end{proof}

\begin{rmk} 
The proof incidentally shows that in the unstably disconnected case, 
$f$ always  has positive Lyapounov exponent with respect to  the Poincar{\'e} metric.
\end{rmk}

If $J$ is disconnected, $f$ is both stably and unstably disconnected. The previous
result suggests the following:

\begin{question}
If $f$ is a polynomial diffeomorphism with disconnected Julia set, is
 it true that almost every connected component of $J$ is a point?
\end{question}

We have no
general answer to this question, except in the case of hyperbolic maps.

\begin{cor}
If $J$ is disconnected and  $f$ is hyperbolic on $J$, 
then almost every component of $J$ is a point. 
\end{cor}

\begin{proof}
if $f$ is hyperbolic on $J$,  it has local product structure. So locally near $p$, $J$ is
homeomorphic to $(J^+\cap W^u_{loc}(p))\times (J^-\cap W^s_{loc}(p))$; 
also the measure $\mu$ is the product of its stable and unstable
conditionals. The previous theorem applied in the stable and unstable 
directions implies the corollary.
\end{proof}

In the next two propositions, we examine the  very special situation
where the submanifolds involved in the decomposition of
$T^-\rest{\bb}$ have bounded degree. In this case we say
$T^-\rest{\bb}$ has {\sl no degree growth}. We do not know whether 
 there are examples exhibiting such a
phenomenon, besides perturbations of polynomials with totally
disconnected Julia sets (horseshoes).
It seems very likely to us that real mappings on the 
boundary of the horseshoe locus indeed satisfy this assumption (see 
\cite[Prop. 3.3]{bs10}).  
We also believe that the condition should hold when
$f$ is hyperbolic and has no attracting orbits.

\begin{prop} \label{prop:nogrowth}
If $T^-\rest{\bb}$ has no degree growth, then
\begin{enumerate} 
\item $f$ is uniformly expanding with respect to the Poincar{\'e} metric
  along the leaves of $T^-\rest{\bb}$.
\item For every unstable manifold $W^u(p)$, $K^+\cap W^u(p)$ is
  totally disconnected. \\
  As a consequence no component of $\Int(K^+)$
  intersects $J^-$. In particular $f$ has no attracting orbits, neither
  attracting rotation sets.
\item $J=J^*\ (:=\supp \mu )$.
\end{enumerate}
\end{prop}

It is proved in \cite[\S 4]{bs8} that 
the first item of the proposition implies ``quasi-expansion''. 
Quasi-expansion means $f$ is uniformly
expanding with respect to a metric related to subsets of the form 
$W^u(p)\cap \set{G^+\leq 1}$ in unstable manifolds; this is
2-dimensional analogue of semi-hyperbolicity in the sense of \cite{cjy}.  
 
\begin{proof} Assume $T^-\rest{\bb}$ has no degree growth, i.e. all
  the disks occurring in the decomposition of $T^-\rest{\bb}$ have
  degree $\leq K$. The set of horizontal subvarieties of degree $\leq
  K$ is compact for the Hausdorff topology. Moreover if a sequence of
  horizontal holomorphic disks $\Delta_n$ converges to $\Delta$ in the
  Hausdorff topology, consider a 
   sequence of parametrizations $\phi_n:\dd\cv \Delta_n$, normalized
  for instance by $\pi\circ\phi_n(0)=0$. The cluster values of the
  sequence $(\phi_n)$
are finite branched coverings $\phi:\dd\cv\Delta$ (use for instance the
  fact that $\pi\circ\phi_n$ is a Blaschke product; see also
  \cite[Lemma 2.2]{cjy}). In particular, consideration of 
the Euler characteristic shows $\Delta$ is a disk.   
Since $\supp T^-= J^-$, $J^-\cap \bb$ is then the  union 
of a family of
 (possibly singular) horizontal disks of degree $\leq K$, namely all cluster 
values of the disks of the decomposition of $T^-\rest{\bb}$.
Moreover, unstable manifolds being disjoint, 
 by the Hurwitz Theorem any two limiting disks $\Delta_1$ and
  $\Delta_2$ are disjoint or equal. 
So  through $p\in J^-$ there is a unique
disk $\Gamma^u(p)$, of degree $\leq K$. \\

As in the proof of the previous theorem, 
by corollary \ref{cor:Poincare}, 
$f^{-1}:\Gamma^u(f(p))\cv f^{-1}\Gamma^u(f(p))\finc \Gamma^u(p)$ is a
 strict contraction for the Poincar{\'e} metric, with factor $\lambda$
 independent of $p$. Moreover $K^+\cap \Gamma^u(p)$ has uniformly
 bounded diameter. Since $f(K^{+,u}(p))=K^{+,u}(f(p))$ it follows that
 every connected component of $K^+\cap \Gamma^u(p)$ is point. It
 follows that no component of $\Int (K^+)$ can intersect $J^-$. The
 corollaries in {\em (2)} are straightforward.\\
 
 For {\em (3)}, let $\Delta$ be a disk in $\bb$
 subordinate to some $\Gamma^u(p)$, such that 
 $\Delta\cap J^+\neq \emptyset$. Assume all forward iterates 
 $f^n(\Delta)$ remain in some vertical sub-bidisk
 $\bb'\subset\bb$. Then $f^n(\Delta)\subset \Gamma^u(f^n(p))\cap \bb'$, which
 contradicts the strict expansion. 
 
 We deduce that if $p\in J=J^+\cap J^-$, the function $G^+$ is not harmonic in any 
  neighborhood of $p$ in $\Gamma^u(p)$. The function 
$G^+$ being continuous, if $\Gamma_n$
 converges to $\Gamma^u(p)$ with multiplicity $k$ in the sense of
  currents ($k\leq K$), 
 $[\Gamma_n]\wedge dd^c G^+$ is positive in the neighborhood of $p$, thus 
$T^-\wedge dd^cG^+>0$ in the neighborhood of $p$ and $p\in J^*$.
 \end{proof}

If moreover $T^+\rest{\bb}$ has no degree growth, we obtain stronger results.

\begin{prop} \label{prop:nogrowth2}
Assume further $T^-\rest{\bb}$ and $T^+\rest{\bb}$ have 
no degree growth. Then: 
\begin{enumerate}
\item The connected components of $J^+\cap \bb$ (resp. $J^-\cap\bb$)
are horizontal disks of bounded degree.
\item If $M$ is any piece of complex curve, not included in 
$J^+$, then $J^+\cap M$ is totally disconnected. The same holds for $J^-$.
\item $\Int(K^+)=\Int(K^-)=\emptyset$.
\item $J=K$ is totally disconnected. 
\end{enumerate}
\end{prop}

\begin{proof} Let $p\in J$. We showed in the preceding
proposition that $K^{+,u}(p)=\set{p}$, so let $\gamma$ be a small
loop in $\Gamma^u(p)$ surrounding $p$ and avoiding $K^+$, so 
that $\Int(\gamma)\cap K^+$ is open and closed in $K^+\cap
\Gamma^u(p)$
(Here $\Int(\gamma)$ is the bounded connected component of $\gamma^c$).

Since $T^+\rest{\bb}$ has no degree growth, $J^+\cap\bb$ is the disjoint union
of a family of vertical disks $\Gamma^s$ with bounded degree. The disks
$\Gamma^s$ move continuously in the Hausdorff topology. From this
we deduce that the union 
$$\bigcup_{q\in \Int{\gamma}} \Gamma^s(q)$$ is both
open and closed in  $J^+\cap\bb$.
As $\gamma$ shrinks to $p$, these graphs converge to $\Gamma^s(p)$. Hence 
$\Gamma^s(p)$ is the connected component of $p$ in $J^+\cap\bb$. \\

Item {\em (2)} is then obvious. Similarly, if $\Int(K^+)\neq\emptyset$, let 
$L$ be a horizontal line intersecting $\Int(K^+)$. Then $\Int(K^+)\cap L$
is a bounded  open set, and $\fr (\Int(K^+)\cap L)\subset J^+\cap L$ which is totally 
disconnected, a contradiction.\\

For {\em (4)}, let $J(p)$ be the connected component of $J$ containing 
$p$. Then $J(p)\subset \Gamma^s(p)\cap\Gamma^u(p)$ which is a finite set. 
We conclude that $J(p)=\set{p}$.
 \end{proof}

\section{Hyperbolic maps with connected $J$}\label{sec:hyp}

In this last section we give an analogue to the familiar result in one dimensional
dynamics, that for hyperbolic maps connectedness is related to the existence 
of attracting orbits. Here the Jacobian determinant tells whether attracting or repelling
orbits must arise. 

\begin{thm}\label{thm:hyperbolic}
Let $f$ be a hyperbolic polynomial diffeomorphism of $\cd$, with 
$\abs{{\jac} (f)}\leq 1$. Assume $J$ is connected. Then $\Int(K^+)$ is
not empty, that is, $f$ has attracting periodic orbits.
\end{thm}

\begin{cor}[\cite{bs7} Corollary A.3] If $f$ is hyperbolic with connected  Julia set, then 
$f$ is not conservative, i.e. $\abs{{\jac} (f)}\neq 1$.
\end{cor}

We actually do believe the analogy with one dimensional dynamics is more
complete.

\begin{question}
If $f$ is a hyperbolic polynomial diffeomorphism with $\abs{{\jac} (f)}\leq 1$ 
and no attracting orbits, does it follow that  $J$ is totally disconnected?
\end{question}

Our opinion is that the answer is yes. More precisely we think that for such a 
$f$, the currents $T^{\pm}\rest{\bb}$ should have no degree growth, in the sense of 
propositions \ref{prop:nogrowth} and \ref{prop:nogrowth2} above.\\

Before embarking to the proof we recall some results of \cite{bs7}. As before, 
let $$V^+=\set{(z,w),\ \abs{z}>\abs{w}, \ \abs{z}>R}.$$ It is proved in \cite{ho1}
that for large enough $R$, and for an appropriate choice 
of $(d^n)^{\rm th}$-root, the sequence $(f^n)^{\unsur{d^n}}$ converges
in $V^+$ to a holomorphic function $\varphi^+$ (the ``B{\"o}ttcher coordinate'') such that
$\varphi^+\circ f = \varphi^d$ and 
\begin{equation}\label{eq:proj}
\varphi^+(z,w) = z+ O(1) \text{ as } V^+\ni z\cv\infty.
\end{equation}
Hence $\varphi^+$ may serve as an invariant first projection at infinity. We saw in section 
\ref{sec:conn} that when $f$ is unstably connected, the unstable leaves are graphs
over the first coordinate near $\fr^v\bb$ for large $R$. By (\ref{eq:proj}),
Rouch{\'e}'s Theorem implies 
these are also graphs for the invariant projection $\varphi^+$ --this
is the exact 
meaning of having no unstable critical points. 

The level lines of $\varphi^+$ define a holomorphic foliation in $V^+$, and
pulling back by $f$, we get an invariant  holomorphic foliation $\mathcal{F}^+$
in $U^+=\cd\setminus K^+$. If $f$ is uniformly hyperbolic and unstably connected,
the leaves of the unstable foliation are transverse to those of $\mathcal{F}^+$
in $V^+$. By using the Lambda lemma, one gets the following result.

\begin{prop}[\cite{bs7} Prop. 2.7] 
If $f$ is unstably connected and hyperbolic, 
the stable lamination $\el^s$ of $J^+$ and 
the holomorphic foliation $\mathcal{F}^+$ fit together into a lamination 
of $J^+\cup U^+$.
\end{prop}

We study the extension of this foliation to $\pd$.

\begin{lem}
$\mathcal{F}^+$ extends as a holomorphic foliation of 
$\overline{V^+}\subset \pd$ by adding the line at infinity as a leaf.
\end{lem}

\begin{proof} The proof is a simple change of coordinates. Let $[Z:W:T]$ 
be the homogeneous coordinates in $\pd$ such that $z=\frac{Z}{T}$ and 
$w=\frac{W}{T}$. In $V^+$ we put $u=\frac{w}{z}=\frac{W}{Z}$ and 
$v=\frac{1}{z}=\frac{T}{Z}$. In coordinates $(u,v)$, $V^+$ becomes 
a bidisk $B_{u,v}=\set{\abs{u}<1,\ \abs{v}<1/R}$, 
with the line at infinity as $(v=0)$.
By (\ref{eq:proj}) we get 
(with an obvious abuse of notation) 
$$\varphi^+(u,v)=\unsur{v}+\delta(u,v),$$
where $\delta$ is a bounded holomorphic function. By Riemann extension, 
$\delta$ extends as a holomorphic function through the line $(v=0)$. \\

Consider the level line $\set{\varphi^+=\unsur{c}}$ for small $c$. This equation rewrites
as $v(1-cv\delta(u,v))=c$. For fixed $u$, the equation has only one solution in $v$ for 
small $c$, depending holomorphically on $c$. This means that in $B_{u,v}$ the level line
$\set{\varphi^+=\unsur{c}}$ is a graph, tending to $(v=0)$ as $c\cv 0$, which is the desired 
result.
\end{proof}

We can now proceed to the proof\footnote{We thank Charles Favre for providing
this elegant argument using Brunella's theorem.} of theorem
\ref{thm:hyperbolic}. Under the 
hypotheses of the theorem, $f$ is unstably connected. 
Assume $\Int(K^+)=\emptyset$. Then $\el^s\cup\mathcal{F}^+$ becomes
 a lamination of $\cd$, holomorphic outside $J^+$. We prove it is globally 
 holomorphic. 

For this, we may adapt a result of Ghys \cite{gh}, or
 use the following direct argument. A theorem of Bowen and Ruelle
 \cite{br} asserts that $J^+=W^s(J)$ has Lebesgue measure zero. Fix a
 flow box for the  lamination $\mathcal{F}^+\cup \el^s$, intersecting
 $J^+$. 
Then for almost every  transverse holomorphic disk $T$, $J^+\cap T$
 has zero measure. We take two such transversals $T_1$ and $T_2$. The
 holonomy map $h:T_1\cv T_2$ is holomorphic in the full measure
 subset   $T_1\setminus J^+$, and globally quasiconformal, so 
it has $L^2_{loc}$ derivatives. By Weyl's lemma, $h$ is holomorphic.

We conclude that $\mathcal{F}^+\cup \el^s$ defines a holomorphic
foliation in $\cd$. This foliation further extends to $\overline{V^+}
\subset \pd$ by adding the line at infinity, by the preceding lemma. 
Pulling back by $f$ (viewed as a birational map on $\pd$), 
allows to extend it to $\pd\setminus I(f)=\pd\setminus [0:1:0]$. We
have thus obtained a singular holomorphic foliation of $\pd$,
preserved by the dynamics. This is a contradiction due to a theorem 
of Brunella \cite{bru}. \hfill $\square$\\

\end{document}